\documentclass[a4paper,12pt]{article}
\usepackage[latin1]{inputenc}
\usepackage{amssymb}
\title{{\bf  Automatic  Ordinals } }
\author{Olivier Finkel$^1$    ~~~~ and  ~~~~  Stevo Todor{\v{c}}evi{\'c}$^{1, 2}$
\\
\\
  $^1${\it Equipe de Logique Math\'ematique }
\\Institut de Math\'ematiques de Jussieu
 \\CNRS and  Universit\'e Paris 7,  France.
\\finkel@math.univ-paris-diderot.fr
\\ stevo@math.univ-paris-diderot.fr
\\
\\ $^2${\it  Department of Mathematics}
\\ University of Toronto, Toronto, Canada
M5S 2E4.
}

\date{}
\begin{document}

\newtheorem{The}{Theorem}[section]
\newtheorem{Pro}[The]{Proposition}
\newtheorem{Deff}[The]{Definition}
\newtheorem{Lem}[The]{Lemma}
\newtheorem{Rem}[The]{Remark}
\newtheorem{Exa}[The]{Example}
\newtheorem{Cor}[The]{Corollary}
\newtheorem{Notation}[The]{Notation}
\newtheorem{Remark}[The]{Remark}
\newcommand{\vp}{\varphi}
\newcommand{\lb}{\linebreak}
\newcommand{\fa}{\forall}
\newcommand{\Ga}{\Gamma}
\newcommand{\Gas}{\Gamma^\star}
\newcommand{\Gao}{\Gamma^\omega}
\newcommand{\Si}{\Sigma}
\newcommand{\Sis}{\Sigma^\star}
\newcommand{\Sio}{\Sigma^\omega}
\newcommand{\ra}{\rightarrow}
\newcommand{\hs}{\hspace{12mm}

\noi}
\newcommand{\lra}{\leftrightarrow}
\newcommand{\la}{language}
\newcommand{\ite}{\item}
\newcommand{\Lp}{L(\varphi)}
\newcommand{\abs}{\{a, b\}^\star}
\newcommand{\abcs}{\{a, b, c \}^\star}
\newcommand{\ol}{ $\omega$-language}
\newcommand{\orl}{ $\omega$-regular language}
\newcommand{\om}{\omega}
\newcommand{\nl}{\newline}
\newcommand{\noi}{\noindent}
\newcommand{\tla}{\twoheadleftarrow}
\newcommand{\de}{deterministic }
\newcommand{\proo}{\noi {\bf Proof.} }
\newcommand {\ep}{\hfill $\square$}
\newcommand{\borapxi}{{\bf\Sigma}^{0}_{\xi}}
\newcommand{\borel}{{\bf\Delta}^{1}_{1}}
\newcommand{\bormpxi}{{\bf\Pi}^{0}_{\xi}}

\maketitle

\begin{abstract}
\noi  
We prove that the injectively $\om$-tree-automatic ordinals are the ordinals smaller than $\om^{\om^\om}$. 
Then we show that the  injectively  $\om^n$-automatic ordinals, where $n\geq 1$ is an integer, are the ordinals smaller than $\om^{\om^n}$. 
This strengthens a recent result 
of Schlicht and Stephan who considered in \cite{SchlichtStephan11}  the subclasses of  {\it finite word} $\om^n$-automatic ordinals. 
As a by-product we obtain  that the hierarchy of injectively  $\om^n$-automatic structures, $n\geq 1$, which was considered in \cite{Fin-Tod2}, is strict. 
\end{abstract}

\noi {\small {\bf Keywords:}    $\om$-tree-automatic structures;  $\om^n$-automatic structures; ordinals.
 }

\section{Introduction}

An automatic structure is a  relational structure whose domain and relations are recognizable by finite automata reading finite words.  Automatic structures
 have very nice decidability and definability properties and have been much studied in the last few years,
see \cite{BlumensathGraedel00,BlumensathGraedel04,KNRS,RubinPhd,RubinBSL}.
 Blumensath considered in \cite{Blumensath99} more powerful kinds of automata. If we replace automata by tree automata (respectively, B\"uchi automata
reading infinite words, Muller or Rabin tree automata reading infinite labelled trees) then we get the notion of  tree-automatic (respectively,
$\om$-automatic, $\om$-tree-automatic) structures.
Notice that  an $\om$-automatic or $\om$-tree-automatic structure may have uncountable cardinality.
All these kinds of automatic structures have the two following fundamental properties.
$(1)$ The class
of automatic (respectively,   tree-automatic, $\om$-automatic, $\om$-tree-automatic) structures is closed under first-order interpretations.
$(2)$ The first-order theory of an automatic (respectively,   tree-automatic, $\om$-automatic, $\om$-tree-automatic) structure is decidable. 

A natural problem is to classify firstly automatic  structures 
using some invariants. For instance Delhomm\'e proved that  the automatic ordinals are the ordinals smaller
than $\om^\om$,  and that the tree-automatic ordinals are the ordinals smaller than $\om^{\om^\om}$, 
see \cite{Delhomme,RubinPhd,RubinBSL}. Kuske proved in \cite{Kuske10} that the $\om$-automatic ordinals are the automatic ordinals, i.e. the ordinals 
smaller than $\om^\om$.  

In the first part of this paper we characterize the (injectively) $\om$-tree automatic ordinals, proving that they are  
the ordinals smaller than $\om^{\om^\om}$. This seems to be the first complete characterization of a class of  $\om$-tree automatic structures. The proof 
uses some  results of Niwinski  \cite{Niwinski91} on the cardinality of regular tree languages, a recent result  of  Barany, Kaiser and Rabinovich \cite{BKR09} on 
(injectively) $\om$-tree automatic structures, the  result of Delhomm\'e on tree automatic ordinals, and some set theory.   

The $\om^n$-automatic structures, presentable by automata reading ordinals words of length $\om^n$, where $n\geq 1$ is an integer, have been recently
investigated in \cite{Fin-Tod2}. 
We show here that the  injectively  $\om^n$-automatic ordinals  are the ordinals smaller than $\om^{\om^n}$. 
This strengthens a recent result 
of Schlicht and Stephan who considered in \cite{SchlichtStephan11}  the subclasses of  {\it finite word} $\om^n$-automatic ordinals. 
As a by-product we obtain  that the hierarchy of injectively  $\om^n$-automatic structures, $n\geq 1$, which was considered in \cite{Fin-Tod2}, is strict. 

The paper is organized as follows. In Section 2 we recall basic notions 
 and some properties of automatic structures. We  characterize the (injectively) $\om$-tree automatic ordinals in Section 3. 
We consider the $\om^n$-automatic ordinals in Section 4. Some concluding remarks are given in Section 5.

\section{Automatic structures}

\noi When $\Si$ is a finite alphabet, a {\it non-empty finite word} over $\Si$ is any
sequence $x=a_1a_2\ldots a_k$, where $a_i\in\Sigma$
for $i=1,\ldots ,k$ , and  $k$ is an integer $\geq 1$. The {\it length}
 of $x$ is $k$.
 The {\it empty word} has no letter and is denoted by $\varepsilon$; its length is $0$.
 For $x=a_1a_2 \ldots a_k$, we write $x(i)=a_i$.
 $\Sis$  is the {\it set of finite words} (including the empty word) over $\Sigma$.

 We introduce now  languages of infinite binary trees whose nodes
are labelled in a finite alphabet $\Si$.
A node of an infinite binary tree is represented by a finite  word over
the alphabet $\{l, r\}$ where $r$ means ``right" and $l$ means ``left". Then an
infinite binary tree whose nodes are labelled  in $\Si$ is identified with a function
$t: \{l, r\}^\star \ra \Si$. The set of  infinite binary trees labelled in $\Si$ will be
denoted $T_\Si^\om$.
A tree language is a subset of $T_\Si^\om$, for some alphabet $\Si$.
(Notice that we shall mainly consider in the sequel {\it infinite} trees so we shall
often use the term tree instead of  {\it infinite} tree).

 We assume the reader has some knowledge about   Muller or  Rabin tree automata and   regular tree languages.
We recall that the classes of  tree languages   accepted by  non-deterministic Muller, Rabin, Street, and parity  tree  automata 
are all the same.  We refer for instance to
\cite{Thomas90,PerrinPin,2001automata} for the definition of these  acceptance conditions.

 Notice that one can  consider a relation 
$R \subseteq T_{\Si_1}^\om \times T_{\Si_2}^\om  \times \ldots \times T_{\Si_k}^\om$, where $\Si_1, \Si_2, \ldots \Si_k$, are  finite alphabets, 
as a tree language over the product alphabet $\Si_1 \times \Si_2 \times \ldots \times \Si_k$.

\hs 

Let now $\mathcal{M}=(M, (R_i^M)_{1\leq i\leq k})$  be  a relational structure,
  where $M$ is the domain,  and for each $i\in [1, k]$ ~ $R_i^M$ is a relation
of finite arity $n_i$ on the domain $M$. The structure is said to be   $\om$-tree-automatic
if there is a presentation of the structure
where the domain and the relations on the domain are accepted by  Muller tree automata,  in the following sense.

\begin{Deff}[see \cite{Blumensath99}]\label{automatic}
Let $\mathcal{M}=(M, (R_i^M)_{1\leq i\leq k})$ be a relational structure, where $n\geq 1$ is an integer,  and each relation $R_i$ is of finite arity $n_i$.
\nl An  $\om$-tree-automatic presentation of the structure $\mathcal{M}$  is formed by   a tuple of  Muller tree  automata
$(\mathcal{A}, \mathcal{A}_=,  (\mathcal{A}_i)_{1\leq i\leq k})$,  and a mapping $h$ from $L(\mathcal{A})$ onto $M$,  such that:
\begin{enumerate}
\ite The automaton $\mathcal{A}_=$ accepts
an equivalence relation $E_\equiv $  on $L(\mathcal{A})$,  and
\ite
For each $i \in [1, k]$, the automaton $\mathcal{A}_i$ accepts an $n_i$-ary relation $R'_i$ on
$L(\mathcal{A})$ such that $E_\equiv$ is compatible with $R'_i$, and
\ite   The mapping $h$ is an isomorphism from the quotient  structure \nl $( L(\mathcal{A}),  (R'_i)_{1 \leq i \leq k} ) / E_\equiv$ onto $\mathcal{M}$.
\end{enumerate}

\noi  The $\om$-tree-automatic presentation is said to be injective if the equivalence relation $E_\equiv $ is just the equality relation on
$L(\mathcal{A})$. In this case  $\mathcal{A}_=$ and $E_\equiv $  can be omitted and
$h$ is simply an isomorphism from $( L(\mathcal{A}),  (R'_i)_{1 \leq i \leq k} )$ onto $\mathcal{M}$.
\noi A relational structure is said to be (injectively) $\om$-tree-automatic if it has an (injective) $\om$-tree-automatic presentation.
\end{Deff}

 We now  recall  a recent decidability result of \cite{BKR09} about injectively $\om$-tree-automatic structures.

 The quantifier $\exists$ is usual in first-order logic. In addition we can consider the infinity quantifier $\exists^\infty$, and the cardinality quantifiers
$\exists^\kappa$ for any cardinal $\kappa$.  The set of first-order formulas is denoted by $\mathrm{FO}$. The set of formulas using the additional quantifiers
$\exists^\infty$, $\exists^{\aleph_0}$ and $\exists^{2^{\aleph_0}}$ is denoted by
 $\mathrm{FO}(\exists^\infty, \exists^{\aleph_0}, \exists^{2^{\aleph_0}})$.
\nl If $\mathcal{L}$ is a set of formulas, the $\mathcal{L}$-theory of a structure $\mathcal{M}$ is the set of sentences (i.e., formulas without free variables)
in $\mathcal{L}$ that hold in $\mathcal{M}$.

If  $\mathcal{M}$ is a relational structure in a signature $\tau$ which contains only relational
symbols and equality, and $\psi$ is a formula in  $\mathrm{FO}(\exists^\infty, \exists^{\aleph_0}, \exists^{2^{\aleph_0}})$,  then
the semantics of the above  quantifiers are defined as follows.

\begin{itemize}
\ite $\mathcal{M} \models \exists^\infty x ~ \psi$  if and only if there are infinitely many $a\in \mathcal{M}$ such that $\mathcal{M} \models  \psi(a)$.
\ite $\mathcal{M} \models \exists^\kappa x ~ \psi$ if and only if the set $\psi^\mathcal{M}=\{a \in \mathcal{M} \mid \mathcal{M} \models  \psi(a) \}$
has cardinality $\kappa$.
\end{itemize}

\begin{The}[\cite{BKR09}]\label{formula}
Let $n\geq 1$ be an integer. If $\mathcal{M}$ is an injectively $\om$-tree-automatic structure 
then its $\mathrm{FO}(\exists^\infty, \exists^{\aleph_0}, \exists^{2^{\aleph_0}})$-theory
 is decidable.  
Moreover assume that  $\mathcal{M}=(M, (R_i^M)_{1\leq i\leq k})$ is a relational structure and  an injective 
$\om$-tree-automatic presentation of the structure $\mathcal{M}$  is formed by   a tuple of  Muller tree  automata
$(\mathcal{A},   (\mathcal{A}_i)_{1\leq i\leq k})$  with an isomorphism $h$ from $L(\mathcal{A})$ onto $M$. If   
 $\psi(x_1, \ldots, x_p)$ is a  formula in $\mathrm{FO}(\exists^\infty, \exists^{\aleph_0}, \exists^{2^{\aleph_0}})$, with free variables 
$x_1, \ldots, x_p$, then the set $h^{-1}[ \{(a_1, \ldots, a_p) \in M^p \mid \mathcal{M} \models  \psi(a_1, \ldots, a_p)\}]$ 
 is a regular subset of  $L(\mathcal{A})^p$ and one can effectively construct a 
tree automaton accepting it.
\end{The}

\section{Injectively $\om$-tree-automatic ordinals}

We  now give a complete characterization of the class of injectively $\om$-tree-automatic ordinals. 

\begin{The}\label{ord-tree}
An ordinal $\alpha$ is an  injectively $\om$-tree-automatic ordinal if and only if it is smaller than the ordinal $\om^{\om^\om}$. 
\end{The}

\proo As usual, an ordinal $\alpha$ is considered as a linear order $(\alpha, <)$. 
We first asssume  that the ordinal $(\alpha, <)$ has an injective $\om$-tree-automatic presentation, which  
  is given by two Muller tree automata $(\mathcal{A}, \mathcal{A}_<)$ 
and an isomorphism $h$ from    $L(\mathcal{A})$ onto $\alpha$.  
 We now distinguish two cases. 

\hs  {\bf First Case.} The ordinal  $\alpha$  is countable. Then   $L(\mathcal{A})$  is countable and then it is proved by Niwinski in \cite{Niwinski91}
 that the infinite trees of the tree language 
$L(\mathcal{A})$  can be represented by {\it finite} binary trees and that the corresponding language of finite trees is regular
 (accepted by a tree automaton reading labelled finite trees). 
Thus the  structure $(\alpha, <)$  is  not only injectively $\om$-tree-automatic but also  injectively tree-automatic, 
the relation $<$ being also accepted by a tree automaton reading 
 {\it finite} binary trees. But Delhomm\'e proved in \cite{Delhomme} that the tree-automatic ordinals 
are the ordinals smaller than the ordinal $\om^{\om^\om}$. Therefore in that case the ordinal 
 $\alpha$  is smaller that the ordinal $\om^{\om^\om}$. 

  On the other hand  it is proved in \cite{CL07} that every 
tree-automatic structure is actually {\it injectively} tree-automatic and hence also {\it injectively}  $\om$-tree-automatic. 
Thus every ordinal smaller that  $\om^{\om^\om}$ is also 
injectively $\om$-tree-automatic. 

Thus the {\it countable} injectively $\om$-tree-automatic ordinals are the ordinals smaller than $\om^{\om^\om}$. 

\hs  {\bf Second Case.} The ordinal  $\alpha$  is uncountable.  Recall that Niwinski proved in \cite{Niwinski91} that a regular language 
of infinite binary trees is either countable or 
has the cardinal $2^{\aleph_0}$ of the   continuum and that one can decide, from a tree automaton $\mathcal{A}$, whether the tree 
language $L(\mathcal{A})$ is countable. 
Since there is an isomorphism $h$ from $L(\mathcal{A})$ onto the uncountable ordinal $\alpha$,  the tree language $L(\mathcal{A})$ and the 
ordinal $\alpha$  have  cardinality $2^{\aleph_0}$. 
In particular,  the ordinal $\alpha$ is greater than or equal to the first uncountable ordinal $\om_1$. But then the initial segment of length
$\om_1$ of $\alpha$ is definable by the formula $\phi(x)= \exists^{\aleph_0} y (y<x) \vee \neg \exists^\infty  y (y<x)$. We can now infer from 
Theorem \ref{formula} that one can effectively 
construct some Muller tree automata  $\mathcal{B}$ and $\mathcal{B}_<$ such that $L(\mathcal{B}) \subseteq L(\mathcal{A})$ and  
$(\mathcal{B}, \mathcal{B}_<)$ 
form an injective $\om$-tree-automatic presentation of the ordinal $(\om_1, <)$. Moreover the formula $\fa x ~\phi(x)$ is satisfied in  the 
structure $(\om_1, <)$ and this can be deduced 
effectively from the two tree automata $\mathcal{B}$ and $\mathcal{B}_<$ of the injective $\om$-tree-automatic  presentation of the structure. 
The ordinal $(\om_1, <)$ being uncountable 
and injectively $\om$-tree-automatic it follows as above that it has the cardinality of the continuum and thus that $2^{\aleph_0}=\aleph_1$, i.e. 
that the continuum hypothesis CH is satisfied. 

We have seen that we can construct 
$(\mathcal{B}, \mathcal{B}_<)$, with an isomorphism $h'$ from   
 $L(\mathcal{B})$ onto $\om_1$,  which 
form an injective $\om$-tree-automatic presentation of the ordinal $(\om_1, <)$. 
Recall that the usual axiomatic system ${\rm ZFC}$ is 
Zermelo-Fraenkel system  ${\rm ZF}$ plus the axiom of choice ${\rm AC}$. 
 A model ({\bf V}, $\in)$ of  the axiomatic system    ${\rm ZFC}$       is a collection  {\bf V} of sets,  equipped with 
the membership relation $\in$, where ``$x \in y$" means that the set $x$ is an element of the set $y$, which satisfies the axioms of  ${\rm ZFC}$.  
We  often say `` the model {\bf V}"
instead of  ``the model  ({\bf V}, $\in)$". The axioms of {\rm ZFC} express some  natural facts that we consider to hold in the universe of sets, and the 
axiomatic system ${\rm ZFC}$ is a commonly accepted framework in which all usual mathematics can be developed. 
We are here supposed to live in a universe {\bf V} of sets which is a model of ${\rm ZFC}$  in which   the continuum hypothesis is  satisfied. 
But we know that, 
using the method of forcing developed by Cohen in 1963 to prove the consistency of the negation of the  continuum hypothesis,  
one can show that there 
exists a forcing extension {\bf V[G]} of {\bf V} which is a model of ${\rm ZFC}$ in which $2^{\aleph_0} > \aleph_1$, see \cite{Jech}.  
Consider now the pair of Muller tree 
automata $(\mathcal{B}, \mathcal{B}_<)$ in this new model {\bf V[G]} of ${\rm ZFC}$. By \cite{Niwinski91} 
the tree language $L(\mathcal{B})$ is   still uncountable (because this can be decided by an algorithm defined in ${\rm ZFC}$ and hence does not depend 
on the ambient model of  ${\rm ZFC}$)
and so it has the cardinality $2^{\aleph_0}$ of the continuum. Moreover the automaton 
 $ \mathcal{B}_<$  defines a linear order $<_0$ on $L(\mathcal{B})$  such that every initial segment is countable (because this is expressed by the  sentence 
$\fa x ~ \phi(x)$ in 
$\mathrm{FO}(\exists^\infty, \exists^{\aleph_0}, \exists^{2^{\aleph_0}})$  which was satisfied in the structure $(\om_1, <)$ in the universe {\bf V}) . 
But one can show that this is in contradiction with the negation $2^{\aleph_0} > \aleph_1$ of the continuum hypothesis. 
Indeed we can construct by transfinite induction  a strictly increasing 
infinite sequence $(x_\alpha)_{\alpha < \om_1}$ in $L(\mathcal{B})$ such that each initial segment 
$I_{\alpha} = \{ x \in L(\mathcal{B}) \mid  x <_0 x_\alpha \}$ is countable. Consider now the union $I=\bigcup_{\alpha < \om_1} I_{\alpha}$. 
The set $I$ forms an initial segment of $L(\mathcal{B})$ for the linear order $<_0$ and it has cardinality $\aleph_1 < 2^{\aleph_0}$. Then there exists 
an element $y\in L(\mathcal{B})$ such that $x_\alpha < y$ for each $\alpha < \aleph_1$. But then the initial segment 
$\{ x \in L(\mathcal{B}) \mid  x <_0 y \}$ is uncountable and this leads  to a contradiction. 
Finally we have proved that  there are no uncountable  injectively $\om$-tree-automatic 
ordinals. 
\ep

\section{The hierarchy of $\om^n$-automatic structures}

Let $\Si$ be a finite  alphabet, and  $\alpha$ be an ordinal; a word of length
 $\alpha$ (or  $\alpha$-word) over the alphabet $\Si$ is an  $\alpha$-sequence $(x(\beta))_{\beta < \alpha}$
(or sequence of length $\alpha$) of letters in $\Si$. The set of $\alpha$-words over the alphabet $\Si$ is denoted by $\Si^\alpha$.
The concatenation of an $\alpha$-word $x=(x(\beta))_{\beta < \alpha}$ and of a $\gamma$-word $y=(y(\beta))_{\beta < \gamma}$ is the 
$(\alpha + \gamma)$-word $z=(z(\beta))_{\beta < \alpha+ \gamma}$  such that $z(\beta)=x(\beta)$ for $\beta < \alpha$ and 
$z(\beta)=y(\beta')$ for $\alpha \leq \beta=\alpha + \beta' < \alpha+ \gamma$; it is denoted $z = x\cdot y$ or simply $z = xy$. 

 We recall now the definition and behaviour of  automata reading words of ordinal length 

\begin{Deff}[\cite{Woj,Woj2,Bedon96}]  An ordinal  Büchi automaton is a  sextuple
 $(\Sigma, Q,$ $q_0, \Delta, \gamma, F)$,  where
$\Si $ is a finite alphabet,
$Q$ is a finite set of states,
$q_0\in Q$ is the   initial state,
$\Delta \subseteq Q \times \Sigma \times Q$ is the transition relation for successor steps, 
$\gamma \subseteq P(Q) \times Q$ is the transition relation for limit steps, and $F   \subseteq Q$ is the set of accepting states.  

 A run of the ordinal   Büchi automaton
 $\mathcal{A}=(\Si, Q, q_0, \Delta , \gamma, F)$  reading a word  $\sigma$ of length
 $\alpha$, is an ($\alpha +1$)-sequence of states  $x$  defined by:
$x(0)=q_0$ and, for  $i<\alpha$, $(x(i), \sigma(i), x(i+1))\in \Delta$ and, for
 $i$  a limit ordinal, $(Inf(x,i),x(i)) \in \gamma$, where $Inf(x, i)$ is the set of states which cofinally appear during the reading
of the  $i$ first letters of  $\sigma$, i.e.
 $$Inf(x, i)=\{ q\in Q \mid  \fa \mu <i, \exists \nu<i  \mbox{ such that }
 \mu<\nu  \mbox{  and  } x(\nu)=q \}$$
\noi A run  $x$  of the automaton  $\mathcal{A}$ over the word
 $\sigma$
of length  $\alpha$ is called successful if $x(\alpha) \in F$. A word $\sigma$
of length  $\alpha$
is accepted  by  $\mathcal{A}$ if there exists a  successful run of $\mathcal{A}$ over $\sigma$. We denote
 $L_\alpha(\mathcal{A})$ the set of words of length  $\alpha$ which are accepted by  $\mathcal{A}$.
 An $\alpha$-language $L$ is a regular  $\alpha$-language if there exists an
 ordinal  Büchi automaton $\mathcal{A}$ such that  $L=L_\alpha(\mathcal{A})$.

An  $\om^n$-automaton is an ordinal  Büchi  automaton reading only words of length $\om^n$  for some integer $n\geq 1$.
\end{Deff}

Recall that we  can obtain  regular  $\om^n$-languages from regular $\om$-languages and regular $\om^{n-1}$-languages 
by the use of the  notion of substitution.

\begin{Pro}[see \cite{Finkel-loc01}]\label{sub} Let $n\geq 2$ be an integer.
An $\om^n$-language $L \subseteq \Si^{\om^n}$  is  regular iff it is obtained from a regular
$\om$-language  $R \subseteq \Ga^\om$
by substituting in every $\om$-word $\sigma \in R$  each letter $a\in \Ga$
by a regular  $\om^{n-1}$-language $L_a  \subseteq \Si^\om$.
\end{Pro}

We can obtain the following  stronger result which will be useful in the sequel. 

\begin{Pro}\label{sub2} Let $n\geq 2$ be an integer.
An $\om^n$-language $L \subseteq \Si^{\om^n}$  is  regular iff it is obtained from a regular
$\om$-language  $R \subseteq \Ga^\om$
by substituting in every $\om$-word $\sigma \in R$   each letter $a\in \Ga$
by a regular  $\om^{n-1}$-language $L_a  \subseteq \Si^\om$, where for all letters $a, b \in \Ga$, $a\neq b$, the languages 
$L_a$ and $L_b$ are disjoint. 
\end{Pro}

\proo The result  follows easily from Proposition \ref{sub} and the two following  facts:
 (1)  The 
class of regular  $\om^n$-languages  is effectively closed under finite union, finite intersection, and complementation, i.e.
we can effectively construct, from two $\om^n$-automata    $\mathcal{A}$ and $\mathcal{B}$, some
$\om^n$-automata
$\mathcal{C}_1$, $\mathcal{C}_2$, and $\mathcal{C}_3$,  such that $L(\mathcal{C}_1)=L(\mathcal{A}) \cup L(\mathcal{B})$,
$L(\mathcal{C}_2)=L(\mathcal{A}) \cap L(\mathcal{B})$,  and $L(\mathcal{C}_3)$ is the complement of $L(\mathcal{A})$.
(2) The class of regular
$\om$-languages $R \subseteq \Ga^\om$ is effectively closed under the substitutions $\Phi: \Ga \ra P(X)$, where $X \subseteq \Ga$ is a finite alphabet and 
$P(X)$ is the powerset of $X$. 
\ep 

\hs We  have defined in \cite{Fin-Tod2} the notion of $\om^n$-automatic presentation of a structure and of $\om^n$-automatic structure 
which is simply obtained by replacing in the above Definition \ref{automatic} the Muller  tree automata by some $\om^n$-automata.

Recently Schlicht and Stephan have independently considered {\it finite word} $\alpha$-automatic structures which are relational structures 
whose domain and relations are accepted by automata reading {\it finite} $\alpha$-words: they define {\it finite } $\alpha$-words as  
words $x$ of length $\alpha$ over an alphabet $\Si$ containing a special symbol $\diamondsuit$ such that all but finitely many letters of $x$ are 
equal to $\diamondsuit$. 
In particular they have proved the following result. 

\begin{The}\label{ord1}
An ordinal  is {\it finite word} $\om^n$-automatic, where $n\geq 1$ is an integer, iff it is smaller than the ordinal $\om^{\om^n}$. 
\end{The}

We now extend this result to the class of injectively $\om^n$-automatic structures considered in \cite{Fin-Tod2}. 

\begin{The}\label{ord2}
An ordinal  is  injectively  $\om^n$-automatic, where $n\geq 1$ is an integer, iff it is smaller than the ordinal $\om^{\om^n}$. 
\end{The}

\proo Notice first that Schlicht and Stephan considered only injective presentations of structures. Thus Theorem \ref{ord1} implies that every ordinal 
$\alpha < \om^{\om^n}$ is  injectively  $\om^n$-automatic. 

On the other hand we proved in \cite{Fin-Tod2} that every injectively $\om^n$-automatic structure is also an injectively $\om$-tree-automatic structure. 
Thus it follows from the above  Theorem \ref{ord-tree} that every  injectively $\om^n$-automatic ordinal is countable. 

But we now prove the following result. 

\begin{Pro}\label{countable}
All  injectively $\om^n$-automatic  countable ordinals  are  {\it finite word} $\om^n$-automatic. 
\end{Pro}

\proo It suffices to show that if $L \subseteq \Si^{\om^n}$ is a   countable regular $\om^n$-language then the 
$\om^n$-words in $L$ can be represented by {\it finite} $\om^n$-words. This is clear for $n=1$ because a countable regular $\om$-language 
$L \subseteq \Si^{\om}$
is of the form $L=\cup_{1\leq j\leq n} U_i\cdot v_i^\om$, where for each integer $i\in [1, n]$  $U_i \subseteq \Sis$ is a finitary language and $v_i$ is a non-empty 
finite word over $\Si$. Then the general result can be proved by an easy induction on $n$, using Proposition \ref{sub2}.
\ep

\hs Finally Theorem \ref{ord2} follows from Theorem \ref{ord1} and Proposition \ref{countable}. 
\ep 

\begin{Cor}
The hierarchy of injectively  $\om^n$-automatic structures, $n\geq 1$, which was considered in \cite{Fin-Tod2}, is strict. 
\end{Cor}

\begin{Rem}
Schlicht and Stephan asked whether a countable  $\alpha$-automatic structure is  {\it finite word} $\alpha$-automatic. The proof of Proposition 
\ref{countable} about ordinals extends easily to the general case of   (injectively) $\om^n$-automatic structures. Thus any countable (injectively) $\om^n$-automatic 
structure is {\it finite word} $\om^n$-automatic. 
\end{Rem}

\section{Concluding remarks}

We have determined the classes of injectively $\om$-tree-automatic ordinals and of injectively $\om^n$-automatic ordinals, where $n\geq 2$ is an ordinal. 
This way we have also proved that the hierarchy of injectively  $\om^n$-automatic structures, $n\geq 1$, which was considered in \cite{Fin-Tod2}, is strict. 

The next step in this research would be to determine the classes of $\om$-tree-automatic and  $\om^n$-automatic ordinals, in the non-injective case. 
Recall that every $\om^n$-automatic structure is a Borel structure, i.e. has a Borel presentation, see \cite{Fin-Tod2},  
and thus it follows from a result of Harrington and Shelah in \cite{HS82}    that there are no uncountable 
$\om^n$-automatic ordinals. We conjecture that a countable $\om^n$-automatic structure is actually injectively  $\om^n$-automatic. 
From this it would follow that the $\om^n$-automatic ordinals are also the ordinals smaller than the ordinal $\om^{\om^n}$. 

\hs {\bf Acknowledgements.}
We  thank  the anonymous referees for  useful comments 
on a preliminary version of this paper.


\begin{thebibliography}{KNRS07}

\bibitem[Bed96]{Bedon96}
N.~Bedon.
\newblock Finite automata and ordinals.
\newblock {\em Theoretical Computer Science}, 156(1--2):119--144, 1996.

\bibitem[BG00]{BlumensathGraedel00}
A.~Blumensath and E.~Gr{\"a}del.
\newblock Automatic structures.
\newblock In {\em Proceedings of 15th IEEE Symposium on Logic in Computer
  Science LICS 2000}, pages 51--62, 2000.

\bibitem[BG04]{BlumensathGraedel04}
A.~Blumensath and E.~Gr{\"a}del.
\newblock Finite presentations of infinite structures: Automata and
  interpretations.
\newblock {\em Theory of Computing Systems}, 37(6):641--674, 2004.

\bibitem[BKR09]{BKR09}
V.~B{\'a}r{\'a}ny, L.~Kaiser, and A.~Rabinovich.
\newblock Cardinality quantifiers in {MLO} over trees.
\newblock In {\em Computer Science Logic, 23rd international Workshop, CSL
  2009, 18th Annual Conference of the EACSL, Coimbra, Portugal, September 7-11,
  2009. Proceedings}, volume 5771 of {\em Lecture Notes in Computer Science},
  pages 117--131. Springer, 2009.

\bibitem[Blu99]{Blumensath99}
A.~Blumensath.
\newblock {\em Automatic Structures}.
\newblock Diploma Thesis, RWTH Aachen, 1999.

\bibitem[CL07]{CL07}
T.~Colcombet and C.~L{\"o}ding.
\newblock Transforming structures by set interpretations.
\newblock {\em Logical Methods in Computer Science}, 3(2), 2007.

\bibitem[Del04]{Delhomme}
C.~Delhomm\'e.
\newblock Automaticit\'e des ordinaux et des graphes homog\`{e}nes.
\newblock {\em Comptes Rendus de L'Acad\'emie des Sciences, Math\'ematiques},
  339(1):5--10, 2004.

\bibitem[Fin01]{Finkel-loc01}
O.~Finkel.
\newblock Locally finite languages.
\newblock {\em Theoretical Computer Science}, 255(1--2):223--261, 2001.

\bibitem[FT12]{Fin-Tod2}
O.~Finkel and S.~Todor{\v{c}}evi{\'c}.
\newblock A hierarchy of tree-automatic structures.
\newblock {\em The Journal of Symbolic Logic}, 77(1):350--368, 2012.

\bibitem[GTW02]{2001automata}
E.~Gr{\"a}del, W.~Thomas, and W.~Wilke, editors.
\newblock {\em Automata, Logics, and Infinite Games: A Guide to Current
  Research [outcome of a Dagstuhl seminar, February 2001]}, volume 2500 of {\em
  Lecture Notes in Computer Science}. Springer, 2002.

\bibitem[HMU01]{HopcroftMotwaniUllman2001}
J.~E. Hopcroft, R.~Motwani, and J.~D. Ullman.
\newblock {\em Introduction to automata theory, languages, and computation}.
\newblock Addison-Wesley Publishing Co., Reading, Mass., 2001.
\newblock Addison-Wesley Series in Computer Science.

\bibitem[HS82]{HS82}
L.~Harrington and S.~Shelah.
\newblock Counting equivalence classes for co-{$\kappa $}-{S}ouslin equivalence
  relations.
\newblock In {\em Logic {C}olloquium '80 ({P}rague, 1980)}, volume 108 of {\em
  Stud. Logic Foundations Math.}, pages 147--152. North-Holland, Amsterdam,
  1982.

\bibitem[Jec02]{Jech}
T.~Jech.
\newblock {\em Set theory, third edition}.
\newblock Springer, 2002.

\bibitem[KL08]{KuskeLohrey}
D.~Kuske and M.~Lohrey.
\newblock First-order and counting theories of omega-automatic structures.
\newblock {\em Journal of Symbolic Logic}, 73(1):129--150, 2008.

\bibitem[KNRS07]{KNRS}
B.~Khoussainov, A.~Nies, S.~Rubin, and F.~Stephan.
\newblock Automatic structures: Richness and limitations.
\newblock {\em Logical Methods in Computer Science}, 3(2):1--18, 2007.

\bibitem[Kus10]{Kuske10}
D.~Kuske.
\newblock Is {R}amsey's {T}heorem omega-automatic?
\newblock In {\em 27th International Symposium on Theoretical Aspects of
  Computer Science, STACS 2010, March 4-6, 2010, Nancy, France}, volume~5 of
  {\em LIPIcs}, pages 537--548. Schloss Dagstuhl - Leibniz-Zentrum fuer
  Informatik, 2010.

\bibitem[Niw91]{Niwinski91}
N.~Niwinski.
\newblock On the cardinality of sets of infinite trees recognizable by finite
  automata.
\newblock In {\em Proceedings of the International Conference MFCS}, pages
  367--376, 1991.

\bibitem[PP04]{PerrinPin}
D.~Perrin and J.-E. Pin.
\newblock {\em Infinite words, automata, semigroups, logic and games}, volume
  141 of {\em Pure and Applied Mathematics}.
\newblock Elsevier, 2004.

\bibitem[Rub04]{RubinPhd}
S.~Rubin.
\newblock {\em Automatic Structures}.
\newblock PhD thesis, University of Auckland, 2004.

\bibitem[Rub08]{RubinBSL}
S.~Rubin.
\newblock Automata presenting structures: A survey of the finite string case.
\newblock {\em Bulletin of Symbolic Logic}, 14(2):169--209, 2008.

\bibitem[SS11]{SchlichtStephan11}
P.~Schlicht and F~Stephan.
\newblock Automata on ordinals and linear orders.
\newblock In {\em Proceedings of the International Conference CiE 2011}, volume
  6735 of {\em Lecture Notes in Computer Science}, pages 252--259. Springer,
  2011.

\bibitem[Tho90]{Thomas90}
W.~Thomas.
\newblock Automata on infinite objects.
\newblock In J.~van Leeuwen, editor, {\em Handbook of Theoretical Computer
  Science}, volume B, Formal models and semantics, pages 135--191. Elsevier,
  1990.

\bibitem[Woj84]{Woj}
J.~Wojciechowski.
\newblock Classes of transfinite sequences accepted by finite automata.
\newblock {\em Fundamenta Informaticae}, 7(2):191--223, 1984.

\bibitem[Woj85]{Woj2}
J.~Wojciechowski.
\newblock Finite automata on transfinite sequences and regular expressions.
\newblock {\em Fundamenta Informaticae}, 8(3--4):379--396, 1985.

\end{thebibliography}
\end{document}